1

# Single Molecule Studies Enabled By Model Based Robust Control Design

Shreyas Bhaban*, *Student Member, IEEE,* Saurav Talukdar*, *Student Member, IEEE,* Mingang Li, Thomas Hays, Peter Seiler, *Member, IEEE* and Murti Salapaka, *Member, IEEE*

*Abstract*—Optical tweezers have enabled important insights into intracellular transport through the investigation of motor proteins, with their ability to manipulate particles at the microscale, affording femto newton force resolution. Its use to realize a constant force clamp has enabled vital insights into the behavior of motor proteins under different load conditions. However, the varying nature of disturbances and the effect of thermal noise pose key challenges to force regulation. Furthermore, often the main aim of many studies is to determine the motion of the motor and the statistics related to the motion, which can be at odds with the force regulation objective. In this article, we propose a mixed objective $H_2/H_\infty$ optimization framework using a model-based design, that achieves the dual goals of force regulation and real time motion estimation with quantifiable guarantees. Here, we minimize the $H_\infty$ norm for the force regulation and error in step estimation while maintaining the $H_2$ norm of the noise on step estimate within user specified bounds. We demonstrate the efficacy of the framework through extensive simulations and an experimental implementation using an optical tweezer setup with live samples of the motor protein 'kinesin'; where regulation of forces below 1 pN with errors below 10% is obtained while simultaneously providing real time estimates of motor motion.

## I. INTRODUCTION

Transport of important cargo inside the cells of eucaryotic species occurs through molecular motor proteins: kineisn, dynein and myosin. These motor proteins enable directed transport of cargo by converting chemical energy to mechanical energy and are central to the regulatory mechanisms that maintain the internal organization of the cell [1], [2]. Structurally, the motor proteins are composed of a cargo binding tail domain, a stalk and two motor heads; where the motor heads bind to the microtubule filaments while transporting the cargo as shown in Fig. 1. The tail domain attaches to the cargo of interest while the motor heads take discrete steps (or 'walk') over the pathways, composed of microtubules, inside the cell [3], thereby enabling directed transport of intracellular cargo. The typical motion of a single motor protein constitutes of series of discrete stepping events [4]. Motor proteins enable crucial intracellular functions, from producing muscle contraction and cellular tension to transport of sub-cellular organelles. Thus, proper functioning of molecular motors is indispensable to maintain a healthy cellular environment. Disruption of transport mechanisms due to impaired motor protein behavior is known to cause a host of maladies, including neurodegenerative diseases such as Huntington's, Parkinson's and Alzheimer's [5], [6], high blood pressure and muscular disorders such as heart disease.

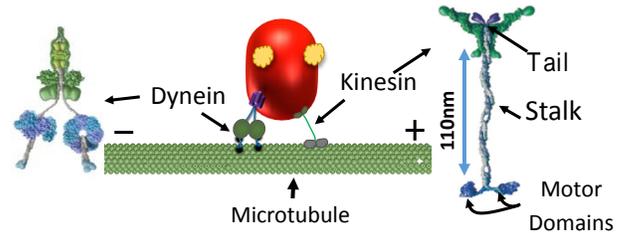

Fig. 1. Schematic showing intracellular cargo carried by motor proteins kinesin and dynein over a section of the microtubule inside the cell. kinesin pulls the cargo towards the 'plus' end of the microtubule while dynein primarily pulls the cargo towards the 'minus' end of the microtubule

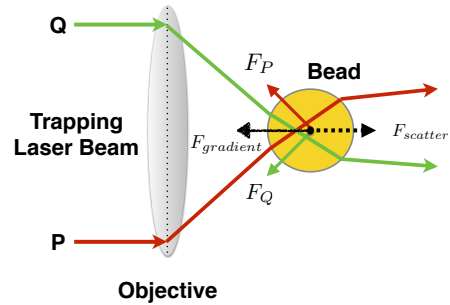

Fig. 2. Bead in an optical trap, as described in [7]. When the trapping laser beam is passed through a high numerical aperture objective, microscopic particles (such as the spherical bead) in the vicinity of the focus of the objective experience two kinds of forces (scattering and gradient forces) due to the momentum transfer onto the particles from the reflected and refracted rays of the laser. $F_P$ and $F_Q$ are the reaction forces felt by the bead due to the momentum transfer from the rays $P$ and $Q$ respectively. The destabilizing *scattering forces* $F_{scatter}$ generated by the laser beam are balanced by the *gradient forces* $F_{gradient}$ resulting from the Gaussian intensity profile of the laser; thereby creating a stable equilibrium point.

Challenges in investigation of motor motion through direct

*Both authors have contributed equally

Shreyas Bhaban (bhaba001@umn.edu, $(612)625 − 3300$) and Murti Salapaka are with the Department of Electrical Engineering, University of Minnesota, Minneapolis, MN, 55455 USA.

Saurav Talukdar are with Department of Mechanical Engineering, University of Minnesota, Minneapolis, MN, 55455 USA

Mingang Li and Thomas Hays are with Department of Genetics Cell Biology and Development, University of Minnesota, Minneapolis, MN, 55455 USA

Peter Seiler is with Department of Aerospace Engineering and Mechanics, University of Minnesota, Minneapolis, MN, 55455 USA

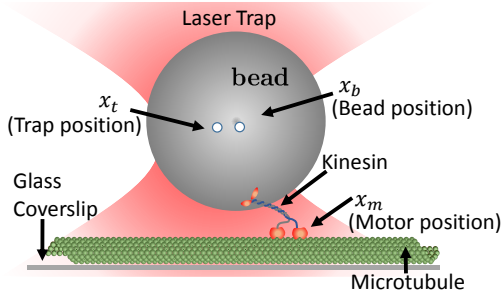

Fig. 3. Investigation of bead attached to motor protein 'kinesin' using an optical trap. $x_t$, $x_b$ and $x_m$ denote the trap location, bead position and motor position respectively.

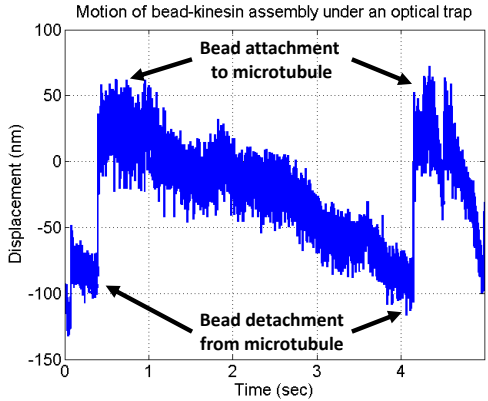

Fig. 4. Bead position data $x_b$ gathered using position sensor, with the bead in an optical trap while simultaneously traversing along the microtubule via the motor protein kinesin attached to it. When the optical trap is held stationary, it exerts more and more force on the bead as the bead is pulled further away from the trap focus. With increasing force on the bead it is likely that the motors pulling the bead detach from the microtubule, at which point the bead is pulled back towards the focus of the trap due to the restoring force of the trap. This figure shows the sawtooth-like patterns [8] characteristic of the bead motion due to the attached molecular motors.

observation (using, for example, confocal or transmission electron microscopy) stem from the extremely small dimensions and step sizes of the motor proteins. For example, kinesin is about 110 $nm$ in length and takes 8 $nm$ sized steps, while exerting forces in the range of femto newtons on the cargo. An efficient method of probing molecular motors is through the use of optical traps [9], [10], which use a laser beam focused through a high numerical aperture (NA) objective to trap micron-sized dielectric beads suspended in a suitable liquid medium, as shown in Fig. 2. When the beam is passed through a high NA objective, the interaction of the bead with the reflected and refracted rays of the laser, results in formation of a stable equilibrium point near the focus of the objective [11]. For small values of bead displacement, $\Delta x$, from the equilibrium point, the trap exerts a restoring force $F_{trap}$ on the bead that is directed towards the trap location. $F_{trap}$ varies linearly with the displacement following the relationship $F_{trap} = K_t \Delta x$; where $K_t$ is the stiffness of the trap. The displacement of the bead can be measured up to a limited range but with sub-nanometer resolution using photo diodes whereas high bandwidth cameras yield larger range with resolution in micrometers [12], [13]. The position of the bead can be controlled by manipulating the position of the trapping laser, using actuators such as galvo mirrors (with bandwidth in Hz) and acousto-optic deflectors (with bandwidth in KHz). Optical traps have enabled force resolution on the order of femto-newtons and position resolution on the order of nanometers, proving successful in the study of a variety of nano-scale systems such as transport inside cells [9], [14], separation of microscopic objects [15], [16], etc.

To facilitate the investigation of motor proteins and detection of the stepping motion (denoted by $x_m$ in Fig. 3) using optical traps, motor proteins are attached to spherical dielectric beads using appropriate biochemistry (see [17]). The beads act as cargo that can be trapped using an optical trapping setup. An *in-vitro* environment mimicking the habitat inside the cell is created in a glass chamber, where microtubule filaments are coated at the base of the glass chamber. Using the principles of optical trapping (as shown in Fig. 2), the cargo is trapped and brought close to the microtubules, to allow for the motor proteins connected to the cargo to attach to the microtubule as shown in Fig. 3. Once attached, the motor protein takes steps on the microtubule, while pulling the bead. An instance of the bead position measured by a photo-diode sensor is shown in Fig. 4. The stepping data collected enables inference of important statistics of molecular motors such as motor detachment rates, stepping rates and reattachment rates, which inform several widely used models of motor proteins. The models have enabled numerous important results on intracellular transport by single and multiple motor proteins [18], [19], [20], [21], [22].

Although measured bead position (such as that shown in Fig. 4) yields stepping statistics of the attached motor proteins, the motor linkage stiffens under increasing load [23], [24] exhibiting nonlinear force extension characteristics. Thus, bead displacement data yields an inaccurate measure of the motor-protein motion, often necessitating corrections (about 15% for forces beyond 1 pN) to compensate for the effect [24]. Furthermore, for bead position measurements such as those shown in Fig. 4, the restoring force on the bead varies as it travels under the influence of the trap. Thus the motor takes every step under changing load conditions. Transport properties of motor proteins are known to be altered when loading conditions change; thus making it difficult to study the effect of load forces on motor protein behavior.

Studies of molecular motors under constant forces using optical tweezers are made possible through the use of *constant force clamps* ([24], [25]). They that are designed to make the optically trapped bead follow the motor protein motion by regulating the separation between bead and trap position, $(x_b - x_t)$ to a desired value $e_d$. allowing for the approximation of motor linkage as a Hookean spring. It is equivalent to maintaining load force on the bead equal to $F_d = k_t e_d$, thus allowing for constant force studies. Furthermore, it enables the measured bead position to correctly reflect the motor motion; from which the motor stepping statistics can be extracted using a variety of offline step detection techniques [26], [27] without the need for ad-hoc corrections [24]. Constant force clamps manipulate the position of the trap $x_t$ using suitable actuators.

The major challenges in the force regulation objective are due to the high bandwidth disturbance caused by motor motion and the large variance of the thermal noise affecting the cargo. State-of-art force clamps [24], [25] yield low bandwidths and are effective for studying motors moving at slow speeds. As alluded to earlier, elucidating the true motion of the motor protein remains the main objective; however, the measured bead motion is corrupted by thermal noise effects and thus presently the measured position of the bead is processed via offline methods [26], [27] to estimate the motor motion. Furthermore, existing methods regulate the load force on the bead but fail to take into account the affect of motor motion on the force experienced by the molecular motor itself. Thus for studies of motor protein based transport of cargo a method that simultaneously provides, accurate loading of the motor protein while it is in motion, and a high bandwidth and high resolution real time estimate of motor motion, will be a significant enabler.

*Our Contribution:* In this article, we employ the mixed objective $H_2/H_\infty$ framework described in [28] to address the the desired diverse set of objectives. Using the framework, it is possible to have regulation of load force as well as estimation of stepping motion as design objectives. Regulation of load force is attained by reducing the $H_\infty$ norm of the transfer function from the disturbance $x_m$ and set point $F_d$ to the regulation error. The effect of motor motion on the force experienced by the motor is minimized by reducing the $H_\infty$ norm of the transfer function from the $x_m$ to error in force experienced by the motor. Accurate estimation of motor motion is attained by minimizing the $H_\infty$ norm of the transfer function from the disturbance $x_m$ to the estimation error, as well as reducing the effects of thermal noise on step estimate by minimizing the $H_2$ norm of the transfer function from thermal noise (white noise) to the motor motion estimate. Through simulations, we demonstrate the ability to regulate forces below 1 pN with errors below 10% on cargoes moving up to speeds of 160 $nm/s$, which is close to the observed kinesin *in-vitro* speeds. Such a capacity will enable constant force studies where low errors in force regulation is of critical importance [21], [20]. Further simulations demonstrate the ability to regulate forces on motors moving as fast as 1000 $nm/s$ with errors below 17%, resulting in significantly improved force clamp capabilities for probing motors with larger step sizes (hence higher speeds [29]) and motors moving at native speeds [30]. Moreover, we demonstrate real time estimation of motor motion for motors travelling at speeds up to 160 nm/s with the root mean squated error in step estimation below 5%; while simultaneously maintaining sub pico-newton forces with regulation error below 10% (as mentioned earlier).

We implement the design on a custom designed optical tweezer setup and demonstrate the performance on live samples of motor protein kinesin. Experimentally, we are able to regulate forces below 1 pN with errors below 8% on motors moving up to 200 nm/s using *in-vitro* motility assays. Without sacrificing force regulation, we demonstrate real time stepping estimation of motor motion at 55 nm/s. This constitutes the first such demonstration, in both simulations and experiments, of maintaining load forces below 1 pN with errors less than 10% while simultaneously estimating motor motion in real time. A preliminary simulation study of the force clamps using robust control framework has appeared in [31] and a preliminary experimental demonstration of the force clamp on a system mimicking the motor motion (without motor proteins) has appeared in [32]. This article builds upon [31], [32] while accommodating the additional objective of reducing the effect of motor motion on the force experienced by the molecular motor. Along with the added objective, the design in this article is able to achieve a notable improvement in the bandwidth of regulation of force on the bead over that reported in [31], [32]. It further presents a detailed computational and experimental implementation of the constant force clamp on live motor protein (kinesin) samples. It thus provides first such demonstration of a model based control design framework on live motor protein samples.

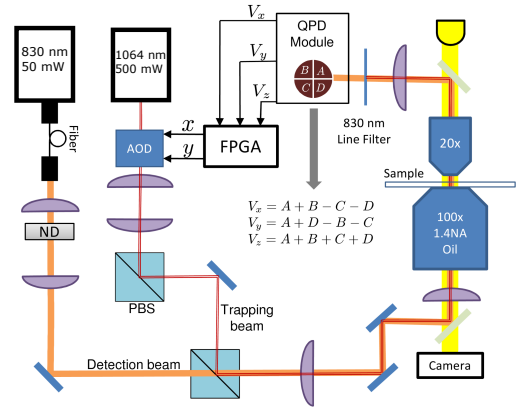

Fig. 5. Schematic for custom designed 'Optical Tweezer' setup

The article is organized as follows. Section II presents the mathematical modeling and characterization of the system. Section III presents the controller design using the mixed objective $H_2/H_\infty$ framework described in [28] for the optical tweezer system. Here, we present simulation results for force regulation and real-time stepping estimation of motor motion in the presence of the thermal noise. Section IV presents experimental implementation of the framework and the results associated with testing the design on live kinesin motor proteins. It is followed by conclusions in section V.

## II. MODELING AND CHARACTERIZATION

### A. Experimental Setup

The custom built optical tweezer forms the experimental setup for realizing optical traps (Fig. 5) where a Nd:YAG trapping laser (CrystaLaser Inc., $\lambda = 1064nm$, $500mW$) is expanded using appropriate optics to fill the back aperture of high numerical apeture (NA) objective (Nikon 100x, 1.4 NA, oil immersion). The optical trap is formed at the focal point where spherical polystyrene beads are trapped. The trapping laser passes through 2-axis acousto-optic-detector (AOD, IntraAction Corp., $DTD - 274HA6$) that precisely steers the beam in $x - y$ plane in response to appropriate commands. To detect the bead position, a second detection laser (Point Source Inc., $iFLEX$ 2000, 50 $mW$, $\lambda = 830nm$) is used to map the image of the bead onto a quadrant photodiode (QPD, Pacific Silicon Sensors, $QP50 - 6SD2$). A neutral density (ND) filter is added in the path of the detection laser to reduce its intensity to ensure that it does not interfere with the trapping phenomenon. The photodiode enables the detection of the bead location by providing three signals $V_x$, $V_y$ and $V_z$ where $V_x$ and $V_y$ represent the photon distribution on the photodiode along $x$ axis and $y$ axis respectively while $V_z$ corresponds to total light intensity. These signals are captured using FPGA based data acquisition system (National Instruments, $PCI7833R$) that generates appropriate commands for the AOD. The controller is implemented using Labview-based NI-FPGA.

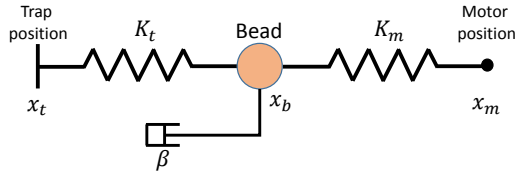

Fig. 6. Block diagram describing the open loop plant

### B. System Model and Instrument Dynamics

Let the locations of the trap centre, cargo and motor head (on the microtubule) be denoted by $x_t$, $x_b$ and $x_m$ respectively. The bead in the optical trap is suspended in an aqueous medium with the coefficient of viscous damping $\beta$. For small magnitudes of displacement $x_b - x_t$ of the optically trapped cargo from the trap centre, the force on the bead due to the optical trap is given by $K_t(x_b - x_t)$ [24]. Similarly, for small extensions $x_m - x_b$ of the motor-stalk, the force on the cargo due to the motor protein is modeled as $K_m(x_m - x_b)$ [24], where $k_m$ is the stiffness of the motor-stalk linkage. Under these conditions, the molecular motor carrying an optically trapped cargo (in this case, a spherical dielectric bead) can be modeled as a spring-mass-damper system with trap stiffness $K_t$, motor stiffness $k_m$ and the coefficient of viscous damping $\beta$, as shown in Fig. 6. The dynamics of the bead are observed to be highly over-damped [7], thus enabling the equation of motion of the bead in the presence of thermal noise $\eta$ to be given by the following Langevian equation,

$$m\ddot{x}_b = -\beta\dot{x}_b + K_m(x_m - x_b) + K_t(x_t - x_b) + \eta. \quad (1)$$

The mass of the bead $m$ is very small (of the order $10^{-17}$ kg), so the Eq. (1) reduces to,

$$\beta\dot{X}_b = K_m(x_m - x_b) + K_t(x_t - x_b) + \eta. \quad (2)$$

Applying Laplace transform to Eq. (2),

$$x_b(s) = G(s)(K_m x_m(s) + K_t x_t(s) + \eta(s)). \quad (3)$$

where $G(s) = \frac{1}{\beta s + K_t + K_m}$ and there is an abuse of notation with $x_m(s)$ used to represent the time and its Laplace transform as well.

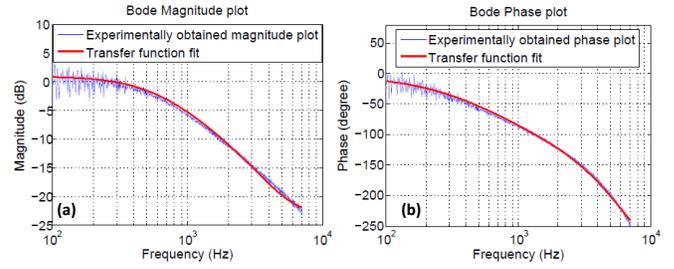

Fig. 7. Experimentally obtained (a) magnitude and (b) frequency plots for transfer function $G_h(s)$ are shown by the 'blue' traces. The 'red' traces indicate a 3-pole 2-zero transfer function fit.

In the absence of the motor (i.e. no motor attached to the bead) and thermal noise, the bead position in open loop follows the equation $x_b(s) = \frac{K_t}{\beta s + K_t}x_t(s)$. The trap position $x_t(s)$ can be manipulated using the command signal to the actuator $u(s)$ through the relation $x_t(s) = A(s)u(s)$, where the transfer function $A(s)$ models the dynamics of the actuator used to control the trap position. Thus, the position of the bead in open loop, in response to the command signal is given by,

$$x_b(s) = K_t G_p(s)A(s)u(s), \quad (4)$$

where $G_p(s) = \frac{1}{\beta s + K_t}$.

### C. Actuator Characterization

One of the primary reasons for the low bandwidth of the state-of-art force clamps is that they ignore the actuator dynamics. The acousto-optic deflector (AOD), which is typically used to manipulate the trap position in optical trapping systems, uses sound waves to create a diffraction grating in a crystal. The nature of the diffraction grating depends upon the frequency of the sound wave used, which is controllable using an input radio frequency (RF) wave. After the sound wave passes through the crystal, it is absorbed by an acoustic absorber. The laser used to create the optical trap is passed through the grating and the first order diffracted spot is utilized to create the optical trap. The spacing of the diffraction grating is altered by changing the frequency of the input RF wave, thereby affecting the position of the trap. Thus the position of the laser spot does not change until the sound wave has crossed the entire laser beam width, making the response time of the AOD dependent upon the speed of the sound wave and the laser beam width. Furthermore, the sound waves suffer from partial reflection from the boundaries of the crystal. Thus, a mixture of waves with different frequencies exists in the crystal after the input frequency is altered. The laser spot settles to the desired position only after the reflected waves completely die out. Therefore, the response time of the AOD is also dependent on the time taken for the reflected waves to die out as well as the absorbing capacity of the acoustic absorber.

Given the inherent physics of the AOD, the dynamics need to be modeled accurately in order to design a force clamp that can accommodate slow as well as fast moving motors. In eq. (4), the transfer function from the command input $u(s)$ to the measured bead position $x_b(s)$ is given by $G_h(s) = K_t G_p(s) A(s)$. We determine the expression for the actuator transfer function $A(s)$ from the experimentally identified expression for $G_h(s)$ using :

$$A(s) = \frac{G_h(s)}{K_t G_p(s)}, \qquad (5)$$

where trap stiffness $K_t$ and damping coefficient $\beta$ are obtained by methods described in [24]. To identify $G_h(s)$ specific to the actuator and control hardware used to perform the experiments in this article, we perform a chirp-wave based system identification where the input signal $u(t)$ is a sinusoidal wave of varying frequencies and output signal is the bead position $x_b(t)$ measured using photo diode sensor. Subsequently, the ratio of amplitudes and the phase difference between the output signal $x_b(j\omega)$ and input signal $u(j\omega)$ is obtained as shown in Fig. 7(a) and Fig. 7(b) respectively. A 3-pole 2-zero transfer function fit yields
$G_h(s) = \frac{0.0396s^2 - 3160s + 1.101 \times 10^8}{1.7 \times 10^{-5}s^3 + 0.9967s^2 + 3.408 \times 10^4 s + 1.086 \times 10^8}$. Fig. 8 shows that the experimentally obtained step response shows a good match with the response predicted by the transfer function $G_h(s)$. Fig. 8 (blue curve) clearly demonstrates that the bead initially traverses in a direction opposite to the commanded direction, which is characteristic of the delays present in the AOD actuation system. The identified transfer function $G_h(s)$ (whose step response is shown in Fig. 8, red curve) contains a right half plane zero, that incorporates effects of delays; and by using (5) yields the actuator transfer function $A(s) = \frac{0.66s^2 - 5.267 \times 10^4 s + 1.835 \times 10^9}{s^2 + 5.51 \times 10^4 s + 1.81 \times 10^9}$. The magnitude and frequency plots of $A(s)$ are shown in Fig. 9 (a) and (b) respectively. Existing force clamps [24], [25] make the simplifying assumption of $A(s) = 1$ which is only appropriate for low frequencies and will lead to diminished performance at higher frequencies.

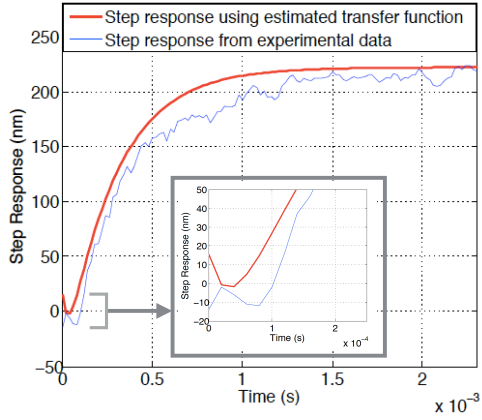

Fig. 8. Validation of the fitted transfer function $G_h(s)$ by comparing step responses. The experimentally obtained step response shows close agreement with that obtained using the fitted transfer function $G_h(s)$. Note the undershoot seen in both the experimental and simulation data (see inset) is characteristic of the delays present in the AOD system.

Another way to identify the actuator transfer function is to model $G_h(s) = e^{-ts}\frac{C}{s+C}$, where $t$ is the actuator delay and the constant $C = \frac{K_t}{\beta}$ incorporates the plant dynamics.

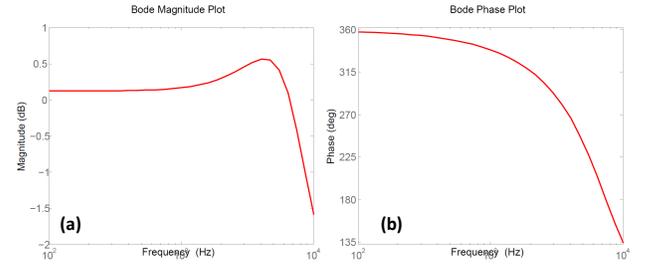

Fig. 9. Magnitude and phase plots for identified transfer function of the actuarot, $A(s)$. Note that at high frequencies, magnitude reduces and phase delay increases, affecting the actuator response for input signals at high frequencies.

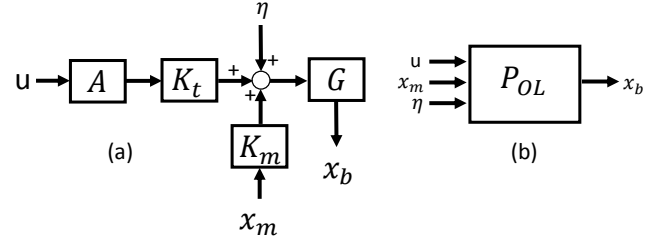

Fig. 10. Block diagram of the open loop plant $P_{OL}$. The inputs to the plant are $u$, $x_m$ and $\eta$ which denote the trap moment command, motor motion and thermal noise respectively. The output $x_b$ denotes the bead position as measured by the photo diode. As mentioned earlier, $G(s) = \frac{1}{\beta s + k_t + k_m}$ and $A(s)$ denotes the transfer function from trap command signal $u(s)$ to the trap position $x_t(s)$, thereby capturing the dynamics of the actuator (AOD).

A second order pade approximation of $e^{-ts}$ yields a transfer function that matches the transfer function $A(s)$ obtained by the fitting performed in Fig. 7. The plant $G(s)$, actuator $A(s)$, trap stiffness $K_t$ and motor linkage stiffness $K_m$, together form the open loop plant $P_{OL}$, as shown in Fig. 10.

In the next section, we present the mixed objective $H_2/H_\infty$ framework utilized in order to design a controller that meets the multiple objectives of maintaining constant force on the bead and estimating the disturbance signal (motor stepping motion).

## III. CONTROLLER SYNTHESIS

A block diagram describing the various components of the optical trapping system together with the feedback controller $K$ is shown in Fig. 11 (a). The controller takes the desired force to be maintained, $F_d$, and the measured bead position, $x_b$, as its inputs and provides the command signal to the AOD $u$ and the estimate $\hat{x}_m$ of the disturbance $x_m$ as its outputs; thus, the controller $K(s)$ is a two-input two-output system as shown in Fig. 12. The performance objectives on the closed loop system are as follows:

1) Good force regulation in the presence of set point changes and disturbances

2) Accurate real time estimation of motor motion

The feedback diagram in Fig. 11 contains additional signals $e_f$ and $\tilde{x}_m$ corresponding to the regulation error and estimation error, respectively. Input-output transfer functions appearing in this loop are denoted via subscripts e.g. $T_{e_f F_d}(s)$

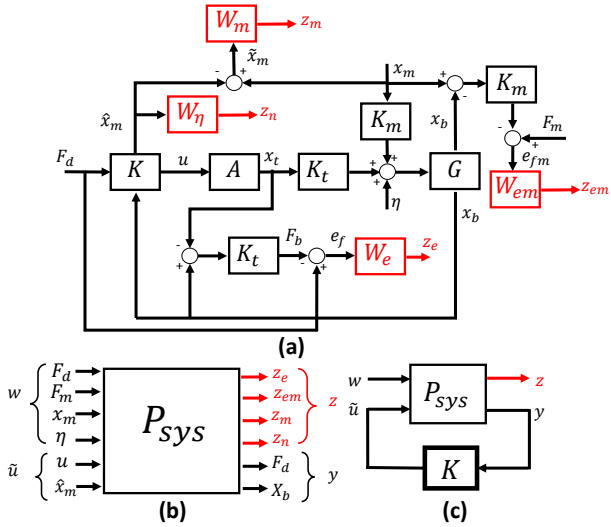

Fig. 11. (a) Block diagram of the plant with frequency dependent weights $W_e, W_{em}, W_m$ and $W_\eta$. (b) $P_{sys}$ denotes the transfer function of the system with inputs $[w; \tilde{u}]$ and outputs $[z; y]$. (c) Closed loop representation of the system with controller $K$.

denotes the transfer function from input $F_d$ to output $e_f$. The desired performance objectives can be expressed in terms of specific transfer functions. The first objective translates to a requirement of small magnitudes for the transfer function $T_{e_f F_d}(s)$ and $T_{e_f X_m}(s)$. Furthermore, the gain of the transfer function $T_{e_f m x_m}(s)$ is maintained to a small value, enabling the regulation of load force on the motor in the presence of motor motion. It also minimizes the effect of motor motion on extensions of the motor linkage, allowing for the the motor linkage to be approximated as a Hookean spring. The second objective translates to a small gain for the transfer function $T_{\tilde{x}_m x_m}(s)$. The presence of thermal noise in the plant output $x_b$ necessitates filtering out of the effect of the white noise from the disturbance estimate $\hat{x}_m$. This translates to reducing the effect of $\eta$ (white noise) on $\hat{x}_m$, that is, minimize the gain of $T_{\hat{x}_m \eta}(s)$ across all frequencies. Note that the disturbance $x_m$ and the noise $\eta$ enter at the same location as shown in Figure 11 (a), so minimizing the $H_2$ norm of $T_{\hat{x}_m \eta}(s)$ is the same as minimizing the $H_2$ norm of $T_{\hat{x}_m x_m}(s)$. Another point to be noted is that, in the optical trapping system, the trap position $x_t$ cannot be directly estimated for fast stepping motors (since $\frac{x_t(s)}{u(s)} = A(s) \neq 1$ at high frequencies); hence, it is not possible to design an effective feedback strategy that uses regulation error $e_f = F_d - k_t(x_b - x_t)$ as an input to the controller. The proposed framework enables desired objectives to be obtained without resorting to the regulation error $e_f$.

### A. Design Objectives

The multi-objective output feedback control paradigm described out in [28] provides an effective controller synthesis approach to meet the desired objectives. The list below enumerates the desired performance requirements :

1) minimize $\|[T_{e_f F_d} \; T_{e_f x_m}]\|_{H_\infty}$

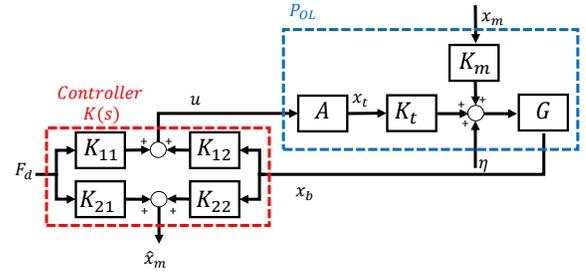

Fig. 12. Block diagram describing the closed loop plant with open loop plant $P_{OL}$ and controller $K(s)$.

2) minimize $\|[T_{e_{fm} x_m}]\|_{H_\infty}$
3) minimize $\|[T_{\tilde{x}_m x_m}]\|_{H_\infty}$
4) minimize $\|[T_{\hat{x}_m x_m}]\|_{H_2}$

We incorporate a constraint on the $H_2$ norm $\|T_{\hat{x}_m x_m}\|_{H_2}$ as $\|T_{\hat{x}_m \eta}\|_{H_2}^2 < \nu$ in the optimization problem, where $\nu$ is specified by the user. Frequency dependent weights $W_e(s)$, $W_{em}(s)$, $W_m(s)$ and $W_\eta(s)$ are introduced to penalize the signals $e_f$, $e_{fm}$, $\tilde{x}_m$ and $\hat{x}_m$ respectively as shown in Figure 11 (a) and ensure feasibility of the optimization problem. Note that, $W_m(s)$ and $W_\eta(s)$ should emphasize non-overlapping frequency regions since the signals $x_m$ and $\eta$ enter the plant at the same location. We define the vector $z := [z_e, z_{em}, z_m, z_n]^T$ as the generalized output vector, $w := [F_d, F_m, x_m, \eta]^T$ as the vector of generalized input, $y := [F_d, x_b]^T$ as the input to the controller, $\tilde{u} := [u, \hat{x}_m]^T$ as the output of the controller $K$. $P_{sys}$ is the open loop plant including the weights, as shown in Figure 11 (b) and is given as, $P_{sys} =$

$$
\begin{bmatrix}
W_e & 0 & -W_e K_t K_m G & -W_e K_t K_m G & W_e K_t A(1-K_t G) \\
0 & W_{em} & -W_{em} K_m (1-K_m G) & W_{em} K_m G & 0 \\
0 & 0 & W_m & 0 & -W_m \\
0 & 0 & 0 & 0 & W_\eta \\
1 & 0 & 0 & 0 & 0 \\
0 & 0 & G K_m & G & G K_t A
\end{bmatrix}.
$$

The associated optimization problem to synthesize a controller $K(s)$ is given below :

$$
\min_{K(s), \gamma} \gamma
$$

$$
\text{such that,} \begin{aligned} \|[T_{z_e f_d} \; T_{z_e x_m}]\|_{H_\infty} &< \gamma \\ \|[T_{z_{em} x_m}]\|_{H_\infty} &< \gamma \\ \|[T_{z_m m}]\|_{H_\infty} &< \gamma \\ \|[T_{z_\eta \eta}]\|_{H_2}^2 &< \nu \end{aligned} \quad (6)
$$

The above optimization problem is a convex optimization problem as shown in [28], the optimal solution to which is obtained by solving the corresponding LMIs. The solution scheme in [28] utilizes a common Lyapunov function for all the objectives, thereby bringing some amount of conservatism in the design process. However, the obtained controller is guaranteed to be stabilizing and also is of the same order as the plant, making it suitable for real time implementation. The associated closed loop system is shown in Fig. 11 (c).

The weight $W_e(s)$ is chosen (Fig. 13(a)) to penalize the regulation error $e_f$, with gains below $-16 \; dB$ for disturbances upto 700 Hz, corresponding to average motor velocity of motors of 5600 $nm/s$. It ensures small values of regulation errors at steady state as well as the frequencies of interest, corresponding to the observed native speeds of molecular motors [30]. $W_{em}(s)$ is chosen (Fig. 13(b)) to minimize the effect of motor disturbances on error $e_{fm}$, by maintaining the gain below $-5\ dB$ for disturbances upto 100 Hz corresponding to an average velocity of 800 $nm/s$, consistent with unloaded speeds of motor protein kinesin [19], [20]. $W_m(s)$ is chosen (Fig. 13(c)) to minimize the effect of disturbances on estimation of motor motion, for motor speeds upto 160 $nm/s$, consistent with observed *in-vitro* speeds using the motility protocol designed in [17] (further details are provided in next section). $W_\eta(s)$ is chosen (Fig. 13(d)) to penalize the higher frequencies ($> 20$ Hz) so that the estimate $\hat{x}_m$ is devoid of high frequency noise. $W_\eta(s)$ is also chosen such that it avoids an overlap of frequencies with the chosen weight $W_m(s)$, while also ensuring the feasibility of the optimization problem. The $2 \times 2$ controller $K(s) = \begin{bmatrix} K_{11}(s) & K_{12}(s) \\ K_{21}(s) & K_{22}(s) \end{bmatrix}$ (shown in the closed loop diagram in Fig. 12) is determined using the given frequency dependent weights and $\nu$ (= 2.23) by solving the corresponding LMI's shown in [28] (using CVX [33]). The optimal value of $\gamma$ is obtained to be 2.63.

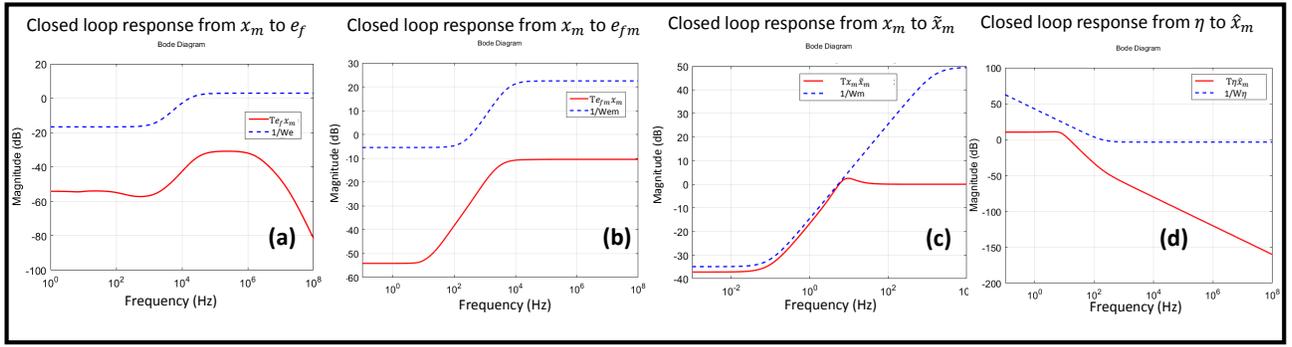

Fig. 13.  Frequency response of (a) $T_{e_f x_m}(s)$ and $W_e^{-1}(s)$,(b) $T_{e_{fm} x_m}(s)$ and $W_{em}^{-1}(s)$, (c) $T_{\hat{x}_m x_m}(s)$ and $W_m^{-1}(s)$, (d) $T_{\hat{x}_m x_m}(s)$ and $W_\eta^{-1}(s)$ using mixed objective $H_2/H_\infty$ method. All three closed loop transfer functions meet the desired frequency dependent performance for the chosen weights $W_e = \frac{0.5s+6000}{s+60}$, $W_{em} = \frac{0.075s+2700}{s+1440}$, $W_m = \frac{0.01s+100}{3s+1.8}$ and $W_n = \frac{1.194\times10^{-3}s}{7.958\times10^{-4}s+1}$.

The weighted closed loop transfer functions of interest, defined in the optimization problem (6) are,

$$
\begin{aligned}
T_{z e_f d} &= W_e(s)(1 - GK_t^2 K_{11}(s)A(s)H(s) \\
&\quad + A(s)K_t K_{11}(s)H(s)), \\
T_{z e x_m} &= W_e(s)K_t G(s)K_m((A(s)K_{12}(s) - 1)H(s)), \\
T_{z_{em} x_m} &= W_{em}(s)K_m G(s)K_t A(s)K_{12}(s)H(s), \\
T_{z_m x_m} &= W_m(s)(1 - G(s)K_m K_{22}(s))H(s),\ \text{and}, \\
T_{z_\eta \eta} &= W_\eta G(s)K_{22}(s)H(s).
\end{aligned}
$$

where $H(s) = \frac{1}{1 - K_t G(s)A(s)K_{12}(s)}$. Note that, $K_{21}(s)$ corresponds to the closed loop transfer function from $F_d$ to $\hat{x}_m$. It is evident that $\|[T_{z e_f d}\ T_{z e x_m}]\|_{H_\infty}$, $\|[T_{z em x_m}]\|_{H_\infty}$, $\|T_{z_m x_m}\|_{H_\infty}$ and $\|T_{z_\eta \eta}\|_{H_2}^2$ do not depend on the controller transfer function $K_{21}(s)$, implying that the estimate $\hat{x}_m$ does not depend on the reference signal $F_d$. Thus, setting $K_{21}(s) = 0$ in the optimal controller obtained by solving the set of LMIs and utilizing controller $\bar{K}(s) = \begin{bmatrix} K_{11}(s) & K_{12}(s) \\ 0 & K_{22}(s) \end{bmatrix}$ does not alter the closed loop performances of the desired channels, as defined in (6). This simplifies the motor motion estimate $\hat{x}_m(s)$ to be $K_{22}(s)X_b$. The frequency response of the inverse of the weights $W_e(s), W_{em}(s), W_m(s), W_\eta(s)$ and the corresponding closed loop transfer functions with the optimal controller in loop are shown in Figure 13. It is clear from the figure that the closed loop transfer functions meet the desired frequency dependent performance specifications.

### B. Simulation results

The performance comparison metric for regulation is the percentage regulation error, $r := \frac{\sigma(F_d - k_t(x_b - x_t))}{F_d} \times 100$, where $\sigma(.)$ denotes standard deviation. The performance metric for estimation error, $e_{rms} := \sqrt{\langle \tilde{x}_m^2 \rangle}$, where $\langle . \rangle$ denotes the expectation operator. Table I lists the regulation and estimation error performance of the synthesized controller for $F_d = 0.95$ pN and $x_m$ being a simulated staircase signal of 8 nm steps with frequencies of $5, 10, 15, 20$ steps/sec corresponding to average motor speeds of $40, 80, 120, 160$ nm/s. It is seen that regulation error is less than 10% and the rms estimation error is less than 6 nm for *in-vitro* motor stepping frequencies upto 160 nm/s. Figure 14 (a) shows an example of the estimate $\hat{x}_m$ of $x_m$ using the noisy bead position $x_b$ when the frequency of the steps is 5 Hz and Figure 14 (b) shows the error in force regulation.

Thus, the synthesized controller achieves the desired objectives of force regulation and real-time motor motion estimation in simulations. In the next section we demonstrate the experimental performance of the controller using live samples of the motor protein kinesin-I.

TABLE I.   REGULATION AND ESTIMATION ERROR FOR VARIOUS MOTOR STEPPING FREQUENCIES

| Stepping Frequency (steps/sec) | Motor Speed (nm/s) | $r$ (%) | $e_{rms}$ (nm) |
|---|---|---|---|
| 5 | 40 | 4.85 | 5.21 |
| 10 | 80 | 7.42 | 4.83 |
| 15 | 120 | 8.06 | 4.78 |
| 20 | 160 | 8.40 | 4.88 |

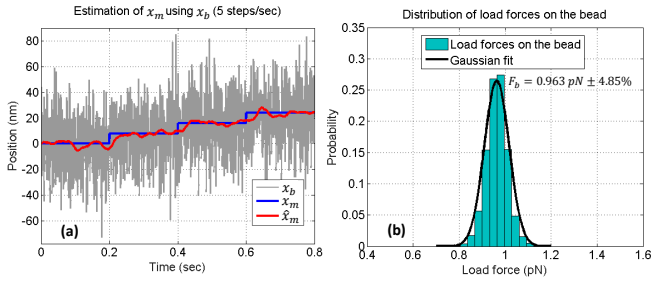

Fig. 14. (a) A realization of $\hat{x}_m$ computed using $x_b$ when simulated $x_m$ is a staircase signal with stepping frequency 5 Hz and (b) the distribution of the associated regulation error. The force on the bead is $F_b = 0.963\ pm \pm 4.85\%$ (1 s.d.)

## IV. EXPERIMENTS WITH LIVE MOTOR PROTEINS

In this section, we present an experimental implementation of the mixed objective $H_2 - H_\infty$ optical force clamp (designed in the previous section) on our optical tweezer setup. The controller is implemented using an National Instruments $PCI7833R$ FPGA. The command signal $u$ to the actuator, which regulates the load force on the bead at the reference value $F_d$ is given by $K_{11}F_d + K_{12}x_b$, where $x_b$ is the measured bead position. The estimated motor motion is obtained by $\hat{x}_m = K_{22}X_b$. The transfer functions $K_{11}$, $K_{12}$ and $K_{22}$ are obtained in the continuous form and have to be discretized in order to be implemented using FPGA hardware. We utilize the bilinear transformation by substituting $s = \frac{2}{T}\frac{z-1}{z+1}$ to obtain the equivalent discrete-time transfer functions $K_{11}$, $K_{12}$ and $K_{22}$, where $T$ is the sampling period of the system (50 KHz in this article). Furthermore, we develop and implement a bead motility assay where carboxylated beads $1\mu m$ in diameter are attached to kinesin-I motor protein using appropriate biochemistry. We then test the performance of the mixed objective clamp on the beads carried by live kinesin-I motors and compare it with the performance of existing state-of-the-art force clamps.

### A. Procedure and demonstration

To investigate motor protein motion using an optical force clamp, the optical tweezer is utilized to trap a polystyrene bead, with the bead being attached to kinesin-I motors using the protocol described in [17]. For the experiment, the trapped bead is brought close to the base of the glass slide where microtubule filaments have been previously coated (See [17] for details) and held stationary at a distance of $\approx 100\ nm$ from it. As the kinesin protein attaches to the microtubule and gradually begins traversing, a unidirectional motion of the bead towards the positively charged end of the microtubule (similar to Fig. 4) is recorded using the photodiode sensor.

The optical force clamp is designed to maintain a constant force of $F_d$ on the bead, which is equivalent to maintaining a constant distance of $d = \frac{F_d}{K_t} = x_b - x_t$ between the locations of the bead and the trap. The clamp is designed to trigger only after the bead has moved more than a distance of $d$ away from the trap focus, that is, after $x_b - x_t \geq d$. When $x_b - x_t < d$, the force clamp is not triggered and the trap position does not change in response to changes in the bead position. An instance of the force clamp in action for $d = 100\ nm$ is demonstrated in Fig. 15 (a), where zones 1, 2 and 3 denote the different stages of bead motion under the optical force clamp. A representation of the bead-motor assembly corresponding to each of the zones is shown in Fig. 15 (b). In zone (1), the bead is trapped and being brought close to the microtubule, awaiting attachment of the motor. No directed motion of the bead is recorded. In zone (2), the motor attaches to the microtubule and pulls the bead, yet the force clamp is not triggered since $x_b - x_t < dx$. In zone (3), $x_b - x_t \geq dx$ and the force clamp is activated, which is seen in the corresponding section in the Fig. 15 (a), where the trap position (red trace) follows the bead position (blue trace) and maintains the distance $x_b - x_t$ close to $d = 100\ nm$. For the trap stiffness $K_t = 0.0095\ pN/nm$, $d = 100\ nm$ translates to a desired force of $F_d = 0.950\ pN$ on the bead.

### B. Results

The mixed objective $H_2 - H_\infty$ force clamp was tested on beads with average velocities ranging from about $59\ nm/s$ to $193\ nm/s$ (see Table II for details). The mixed objective force clamp demonstrated significantly improved regulation performance, with the error in force resolution ranging from $4.67\%$ to $7.97\%$. An instance of the force clamp operating on a bead travelling at $107\ nm/s$ is shown in Fig. 16 (a), where an error of in force regulation of $6.075\%$ is obtained (as shown in Fig. 16 (b)). We compared the performance of the existing state-of-art [24] under similar conditions of load force and velocities, where higher errors in force regulation were seen (as shown in TABLE III). An instance of the traditional force clamp operating on a bead travelling at $149\ nm/s$ is shown in Fig. 17 (a), where an error of $17.68\%$ is obtained (Fig. 17 (b)).It is evident that the mixed objective approach provides an advantage over the existing approaches and is capable of regulating sub $pN$ forces with errors as low as $4.67\%$.

We further utilized the mixed objective force clamp to provide with a real time estimate $\hat{x}_m$ of the motor motion $x_m$. The estimate of the motor motion for an average velocity of 55 nm/s ($\approx$ stepping frequency of 7 steps/sec) is shown by red curve in Figure 18. Since the actual stepping signal is not known, an estimate obtained using state-of-art offline step estimation algorithm [27] (blue curve) is shown for comparison. It is seen that the high frequency noise is filtered out, providing an estimate $\hat{x}_m$ comparable with state-of-art offline step estimation techniques. Moreover, the estimated signal $\hat{x}_m$ can provide starting points or edges that can be used by machine learning or dynamic programming based algorithms to generate statistics in real time or with an acceptable amount of delay. Thus, the force clamp presented in this article facilitates investigation and real time estimation of motor motion under sub $pN$ load forces, which to the best of our knowledge has not been reported in the existing studies.

TABLE II. $H_2 - H_\infty$ FORCE CLAMP REGULATION PERFORMANCE

| Motor Speed in ($nm/s$) | Error in force regulation (%) |
|---|---|
| 59.17 | 4.67 |
| 68.11 | 5.069 |
| 93.2 | 7.07 |
| 107.53 | 6.075 |
| 152.01 | 5.073 |
| 193.3 | 7.97 |

TABLE III. PERFORMANCE USING TRADITIONAL FORCE CLAMP [24]

| Motor Speed in ($nm/s$) | Error in force regulation (%) |
|---|---|
| 45.86 | 17.806 |
| 149.64 | 17.68 |

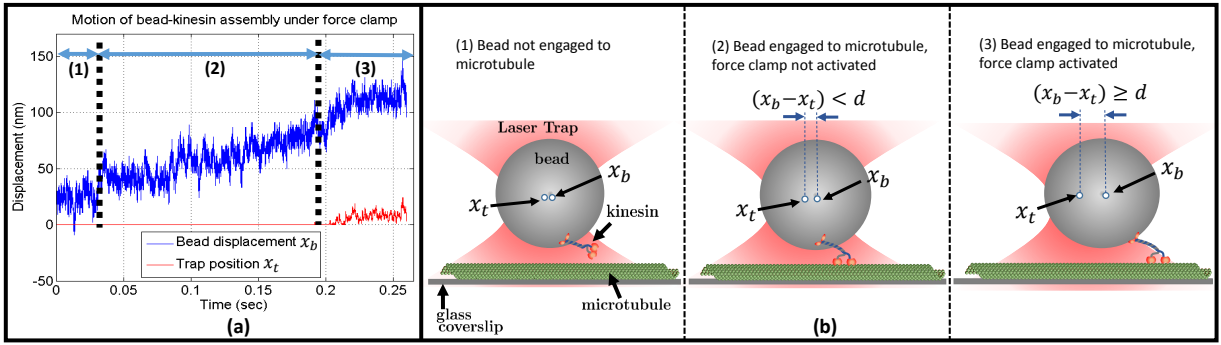

Fig. 15. (a) Sample trace of bead ($x_b$) and trap position ($x_t$) showing optical force clamp in action, divided into three zones. (b) Diagrammatic representation of the behavior of bead-kinesin ensemble in each of the three zones. In zone (1), the bead attached to kinesin is brought close to the glass surface and awaits the attachment of the kinesin to the microtubule filaments coated onto it. No motion of the bead is recorded as the kinesin does not attach to the microtubule. Once the attachment occurs as shown in zone (2), the motor pulls the bead and traverses along the microtubule in a directed manner. However, since the displacement of the bead position from the trap position is less than $d = 100\ nm$, the force clamp is not triggered and the trap position $x_t$ remains unchanged. Once $x_b - x_t > 100\ nm$ as seen in zone (3), the force clamp is triggered and the trap position follows the bead position while maintaining $x_b - x_t$ close to $100\ nm$.

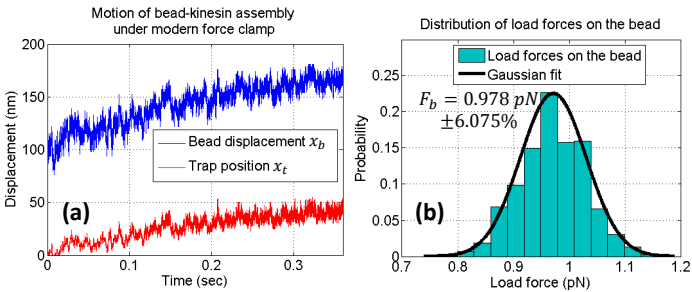

Fig. 16. (a) Sample trace of bead and trap position showing the modern $H_2/H_\infty$ force clamp in action (b) Distribution of the error in force regulation.

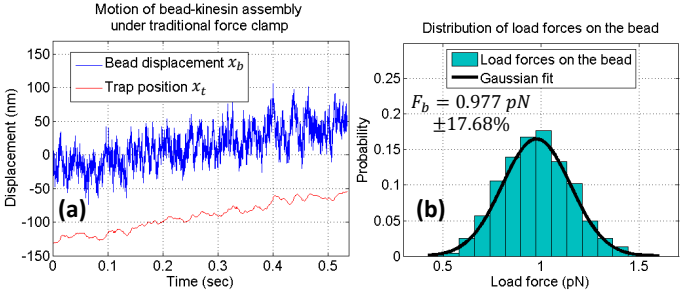

Fig. 17. (a) Sample trace of bead and trap position showing the traditional force clamp [24] in action (b) Distribution of the error in force regulation.

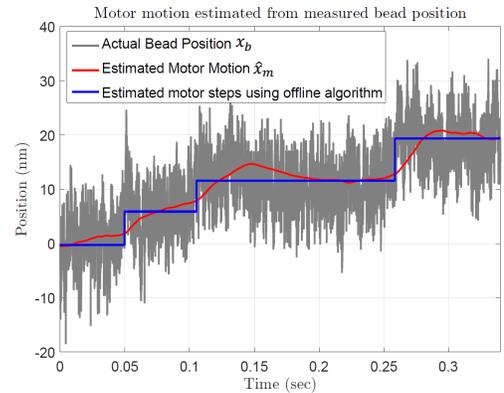

Fig. 18. Real time estimation of stepping motion, with average bead velocity of 55 nm/s

## V. CONCLUSION

This article presents a detailed design for an optical force clamp using a mixed objective $H_2/H_\infty$ framework, along with verifications using simulations as well as experiments using live motor proteins; with extensive applications to investigation of bio-molecules and processive motor proteins. The model-based approach using a systematic design methodology results in notable improvement in force regulation and bandwidth of operation over existing designs [25], [24]. A crucial enabler of the design is a precise understanding and modeling of hitherto un-modeled dynamics in the existing instruments using right half plane zeros; identified by a data-driven system identification approach. We verify the design through extensive numerical simulations, providing improved force regulation with guaranteed performances based on user-defined weight functions. The improved design can enable bio-molecular studies with higher resolution through the ability to investigate systems at finer intervals of forces. It will allow researchers to discover events that were not possible due to limitations of current designs, leading to possible refinement of bio-molecular models. The design allows for examining motor proteins at higher velocities close to their native speeds [30] and proteins with larger step sizes [29] with satisfactory force regulation. The design using modern control framework adds to the robustness of the system to parametric uncertainties. We further implemented the design on an optical tweezer setup and tested its efficacy using live samples of the motor protein 'kinesin', demonstrating the ability to regulate sub-pico newton forces with errors below 10%; individual attempts are also shown to have errors below 5%.

The article also proposes a scheme for real time stepping estimation of molecular motors, which to the best of our knowledge is absent from current literature. By virtue of the mixed objective $H_2/H_\infty$ design, we are able to reduce the effects of thermal noise on stepping estimate without compromising on the force regulation. We demonstrate the efficacy of the method using simulations as well as experiments involving live 'kinesin' samples. Unlike some methods that reduce the effects thermal noise on stepping signal through stiffening the system [34], [24], our method preserves the dynamics of the system and thus avoids impacting the system being examined. The real time step estimation methodology can prove to be of importance to experimental bio-physicists, equipping them with the ability to detect and identify motor protein motility in real time. To the best of the author's knowledge, our method is the first of its kind that demonstrates the force regulation with errors as low as $4.67\%$ on samples of live motor protein kinesin, while simultaneously providing real time estimates of the motor stepping motion. Note that the step estimation bandwidth for our method is limited to slow moving motors due to the trade-off in the $H_2/H_\infty$ cost over the identical channel. It is useful for motor proteins operating in the high load force and low ATP regimes [24], [35]. Future work entails exploiting the knowledge of the noise properties and designing the corresponding weight function ($W_\eta$ in this article) more efficiently. Nevertheless, our proposed scheme enables improved investigative abilities in molecular biology using system theory tools and provides improvement in force clamping abilities over current methods, with the added contribution of estimation of motor motion in real time.

## ACKNOWLEDGMENT

The authors would like to thank Subhrajit Roychowdhury, GE Global Research, for help with the optical trapping experiments and Tanuj Agrawal, Cymer for sharing his insights on optical tweezers. The discussions with Biswaranjan Mohanty and Sivaraman Rajaganapathy from the University of Minnesota have been insightful. The research is supported by National Science Foundation (http://www.nsf.gov/) grant number CNS 1544721.